\documentclass[a4paper,12pt]{article}
\usepackage{amsmath,amsthm,amssymb}
\usepackage[mathscr]{eucal}
\usepackage{enumerate}
\theoremstyle{plain}
\newtheorem{proposition}{Proposition}[section]
\newtheorem{theorem}{Theorem}[section]

\theoremstyle{definition}

\theoremstyle{remark}

\newtheorem{remark}{Remark}[section]
\newcommand{\dom}{\mathop{\rm dom}\nolimits}
\newcommand{\ran}{\mathop{\rm ran}\nolimits}
\newcommand{\rank}{\mathop{\rm rank}\nolimits}

\begin{document}
\author{G.Y. Tsyaputa }
\title{Green's relations on the deformed transformation semigroups}
\date{}
\maketitle
\begin{abstract}
Green's relations on the deformed finite inverse symmetric
semigroup $\mathcal{IS}_n$
    and the deformed finite symmetric semigroup $\mathcal{T}_n$ are described.

\end{abstract}

\section{Introduction} Let $X$ and $Y$ be nonempty sets,
$S$ a set of maps from $X$ to $Y$. Let $\alpha : Y\rightarrow X$
be a fixed map. We define the multiplication of the elements from
$S$ by $\phi\circ \psi = \phi\alpha\psi$ (the compositions of the
maps is from left to right). The action defined above is
associative. E.S.Ljapin (\cite{Ljap},~p.~393) formulated the
problem of investigation of the properties of this semigroup
depending on the restrictions to set $S$ and map $\alpha$.

 Magill \cite{Mag} studied this problem in the case of topological spaces and
 continuous maps. Under the assumption that $\alpha$ be onto he described the automorphisms
 of such semigroups and determined their isomorphism criterion.

Later Sullivan \cite{Sul} proved, if $|Y|\leq|X|$ then Ljapin's
semigroup is embedding into transformation semigroup on the set
$X\cup \{a\}$, $a\notin X$.

An important case is when $X=Y$, $T_X$ is a transformation
semigroup on the set $X$, $\alpha\in T_X$. Symons \cite{Sym}
stated the isomorphism criterion for such semigroups and
investigated the properties of their automorphisms.

The latter problem may be generalized to arbitrary semigroup $S:$
for a fixed $a\in S$ the action $*_a$ is defined by $x*_a y=xay$.
 The obtained semigroup is denoted by $(S,*_a)$ and action $*_a$ is
called the multiplication deformed by element $a$ (or just the
deformed multiplication).

In~\cite{Tsyap} pairwise nonisomorphic semigroups received from
finite symmetric semigroup $\mathcal{T}_n$  and finite inverse
symmetric semigroup $\mathcal{IS}_n$ are classified. In particular
there holds
\begin{theorem}
\label{izo} Semigroups $(\mathcal{IS}_n ,*_a)$ and
$(\mathcal{IS}_n ,*_b)$ are isomorphic if and only if $\rank
(a)=\rank (b)$.
\end{theorem}

In this article Green's relations on $(\mathcal{IS}_n ,*_a)$ and
$(\mathcal{T}_n ,*_a)$  are described.

Recall that $\mathcal{L} -$relation is defined by $a\mathcal{L} b
\Longleftrightarrow S^1 a =S^1 b$. Similarly $\mathcal{R}
-$relation is defined by $a\mathcal{R} b \Longleftrightarrow aS^1
=bS^1$, and $a\mathcal{J} b \Longleftrightarrow S^1 aS^1 =S^1
bS^1$. The following Green's relations are derivative:
$\mathcal{H} = \mathcal{L} \cap\mathcal{R}$, $\mathcal{D} =
\mathcal{L} \vee\mathcal{R}$.

Since on finite semigroups $\mathcal{D}=\mathcal{J}$
\cite[prop.~2.3]{Lal}, $\mathcal{J} -$ relation is not considered
further below. Denote $L_a (R_a ,H_a ,D_a)$ the class of the
corresponding relation containing $a$.
\begin{remark}
\label{umovy}It is obvious, the elements $a,\:b$ belong to the
same $\mathcal{L} -$class (resp. $\mathcal{R} -$class) if and only
if there exist such $u,v$ from $S$, that $a=ub$ and $b=va$ (resp.
$a=bu$ and $b=av$).
\end{remark}
Semigroup $(S,\circ)$ with operation $a\circ b=b\cdot a$ for any
$a,b$ from $S$ is called dual to semigroup $(S,\cdot)$. If
semigroups $(S,\cdot)$ and $(S,\circ)$ are isomorphic than
$(S,\cdot)$ is called self-dual.

We follow terminology and notation as in~\cite{Art}.
\section{Green's relation on $(\mathcal{IS}_n ,*_a)$.}
Let $x\in\mathcal{IS}_n.$ Denote $x^{-1}$ an inverse element to
$x$. If $x$ is represented as a partial permutation
\[x= \left( \begin{array}{lccccr} i_{1} &\ldots& i_{k}
& i_{k+1} & \ldots & i_{n}\\
 j_{1} & \ldots & j_{k} & \varnothing & \ldots & \varnothing \end{array}
 \right) \]
then $x^{-1}= \left( \begin{array}{lccccr} j_{1} &\ldots& j_{k}
& j_{k+1} & \ldots & j_{n}\\
 i_{1} & \ldots & i_{k} & \varnothing & \ldots & \varnothing \end{array}
 \right)$.
For a partial transformation $x\in\mathcal{IS}_n$ denote by $\dom
(x)$ the domain of $x$, and $\ran (x)$ the range of $x$. The value
$|\ran (x)|$ is called the rank of $x$ and is denoted by $\rank
(x)$.

\begin{theorem}
\label{dvoist} Semigroup $(\mathcal{IS}_n ,*_a)$ is self-dual.
\end{theorem}
\begin{proof}
Denote $(\mathcal{IS}_n ,\circ_a) $ the semigroup which is dual to
$(\mathcal{IS}_n ,*_a)$. We show that $(\mathcal{IS}_n ,*_a)$ and
$(\mathcal{IS}_n ,\circ_a) $ are isomorphic to the semigroup
$(\mathcal{IS}_n ,*_{a^{-1}})$.

Semigroups $(\mathcal{IS}_n ,*_a)$ and $(\mathcal{IS}_n
,*_{a^{-1}})$ are isomorphic by the theorem~\ref{izo}. On the
other hand, one can easily check the map
\begin{displaymath}f:(\mathcal{IS}_n,*_{a^{-1}})\Rightarrow (\mathcal{IS}_n
,\circ_a), \;x\mapsto x^{-1}
\end{displaymath} is an isomorphism.
\end{proof}
\begin{proposition}
\label{treba} For any $x,y$ from $(\mathcal{IS}_n ,*_a)$ such
element $u$ from $(\mathcal{IS}_n ,*_a) $ that $x=y*_a u$, exists
if and only if the following two conditions hold:
\begin{displaymath} \dom (x)\subseteq\dom (y);\qquad (1)\qquad\qquad
  \ran (y)\subseteq\dom (a).\qquad (2)
  \end{displaymath}
\end{proposition}
\begin{proof}
Let $x=y*_a u$ for $u$ from $(\mathcal{IS}_n ,*_a).$ Then it is
clear that $\rank (x)\leq\rank (a)$. If $i\in\dom (x)$ then
$x(i)=y*_a u(i)=u(a(y(i)))$, so $i\in\dom (y)$. If $y(i)\notin\dom
(a)$ then $i$ does not belong to the domain of $x$, so $\ran
(y)\subseteq\dom (a)$.

Conversely, let conditions (1)-(2) hold. Consider a partial
permutation, $u$, defined only on the elements from $\ran (a)$,
moreover, $u(i)=x(y^{-1}(a^{-1}(i)))=a^{-1}y^{-1}x(i),\;i\in\ran
(a)$. It is obvious that $y*_a u=yau=x$.
 \end{proof}
\begin{theorem}
\label{klasy} Let $x\in (\mathcal{IS}_n ,*_a)$.
\begin{enumerate}[1)]
    \item\label{1)} If $\ran (x)\subseteq\dom (a)$ then
    \begin{displaymath}R_x = \{y\:|\dom (y)=\dom (x),\: \ran (y)\subseteq\dom
    (a)\};
    \end{displaymath}

    otherwise $R_x =\{x\}$.
    \item\label{2)} If $\dom (x)\subseteq\ran (a)$ then
    \begin{displaymath}L_x = \{y\:|\ran (y)=\ran (x), \:\dom (y)\subseteq\ran
    (a)\};
    \end{displaymath}

    otherwise $L_x =\{x\}$.
    \item\label{3)} If $\ran (x)\subseteq\dom (a)$ and $\dom (x)\subseteq\ran (a)$ then
    \begin{displaymath}H_x =\{y\: |\dom (y)=\dom (x), \: \ran
    (y)=\ran(x)\};
    \end{displaymath}
    otherwise $H_x =\{x\}$.
    \item\label{4)} If $\ran (x)\subseteq\dom (a)$ and $\dom (x)\nsubseteq\ran
    (a)$ then $D_x =R_x$;

   if $\ran (x)\nsubseteq\dom (a)$ and $\dom (x)\subseteq\ran (a)$ then $D_x
   =L_x$;

   if $\ran (x)\subseteq\dom (a)$ and $\dom (x)\subseteq\ran (a)$ then
   \begin{displaymath}D_x =\{y\: |\dom y\subseteq\ran (a),\:
    \ran (y)\subseteq\dom (a)\};
    \end{displaymath}

   otherwise $D_x =\{x\}$.

   \end{enumerate}
In particular if $\rank (a)\leq 1$ then all Green's relations
classes on semigroup $(\mathcal{IS}_n ,*_a)$ are one element.
\end{theorem}
\begin{proof} 1) By remark ~\ref{umovy} for $x$ and $y$ from semigroup $(\mathcal{IS}_n ,*_a)$
belong to the same $\mathcal{R} -$class there should exist such
$u$ and $v$ from $(\mathcal{IS}_n ,*_a) $ that
\begin{gather}
x=y*_a u, \tag{3}\label{$i$} \\ y=x*_a v. \tag{4}\label{$ii$}
 \end{gather}
By Lemma~\ref{treba} equality~(\ref{$i$}) holds if and only if
\begin{displaymath}\dom (x)\subseteq\dom (y);\quad (5)\qquad\qquad \ran
(y)\subseteq\dom (a),\quad (6)
\end{displaymath} and (\ref{$ii$})
holds if and only if
\begin{displaymath}
\dom (y)\subseteq\dom (x); \quad(7)\qquad\qquad \ran
(x)\subseteq\dom (a).\quad (8)
\end{displaymath}
Condition $(8)$ implies if $\ran (x)\nsubseteq\dom (a)$ then $R_x
=\{x\}$. Let $\ran (x)\subseteq\dom (a)$. Following conditions
$(5)$ and $(7)$, one gets  for all $y\in R_x$ there holds an
equality $\dom (x)=\dom (y)$.  Finally condition $(6)$ implies
$\ran (y)\subseteq\dom (a)$. Besides conditions $(5)-(8)$ are
sufficient so statement $1)$ is proved.

Obviously, if $\rank (a)\leq 1$ then $|R_x| =1$.

2) By Theorem~\ref{dvoist} $L_x -$ class description is obtained
from  $R_x -$ class description by the interchanges of domains and
ranges.

3) Statement about $H_x -$ class follows from 1) and 2) statements
of the theorem and $\mathcal{H} -$ relation definition.

4) $x\mathcal{D}y$ if and only if there exists such $z$ from
$(\mathcal{IS}_n ,*_a)$ that  $x\mathcal{L} z$ \: and
$z\mathcal{R} y$. Consider the possible cases.
\begin{enumerate}[a)]
    \item $|L_x|=1$ and $|R_x|>1$. As mentioned above, this case holds if and only
    if $\ran (x)\subseteq\dom (a)$ and $\dom (x)\nsubseteq\ran
    (a)$;
    then $D_x =R_x$.
    \item $|L_x|>1$ and $|R_x|=1$. Then $\ran (x)\nsubseteq\dom (a)$ and $\dom (x)\subseteq\ran
    (a)$,
    hence $D_x =L_x$;
    \item $|L_x|>1$ and $|R_x|>1$. In this case $\ran (x)\subseteq\dom (a)$ and $\dom (x)\subseteq\ran
    (a)$; so $y\mathcal{D} x$ if and only if $\ran (y)\subseteq\dom (a)$ and $\dom (y)\subseteq\ran
    (a)$.
    \item Now if $|L_x|=1$ and $|R_x|=1$ then it is obvious that $D_x
    =\{x\}$.
\end{enumerate}
\end{proof}
\begin{proposition}
Let $p=\rank (a)$, $p>1$.  Then in semigroup $(\mathcal{IS}_n
,*_a)$

the number of one element $\mathcal{R} -$classes ($\mathcal{L}
-$classes) equals
\begin{displaymath}
\sum_{k=0}^{n}\sum_{m=1}^{k}\binom{n-p}{m}\binom{p}{k-m}\binom{n}{k}k!\,.
\end{displaymath}
The number of multi-element $\mathcal{R} -$classes ($\mathcal{L}
-$classes) equals
$\displaystyle{\overset{p}{\underset{k=1}{\sum}}}\binom{n}{k}$.
The cardinality of a multi-element class is
\begin{displaymath}[p]_k=p(p-1)\cdots (p-k+1), \: \text{where} \; 1\leq k\leq p,
\end{displaymath} moreover the number of $\mathcal{R}
-$classes ($\mathcal{L} -$classes) of the cardinality $[p]_k$
equals $\binom{n}{k}$.
\end{proposition}
\begin{proof}
By Theorem~\ref{dvoist} semigroup $(\mathcal{IS}_n ,*_a)$ is
self-dual so the number and the cardinalities of $\mathcal{R}
-$classes as well as of $\mathcal{L} -$classes are the same. Hence
it is enough to calculate the number and the cardinalities of
$\mathcal{R} -$classes.

Let $\rank (x) =k$.

The proof of the statement $1)$ of the Theorem~\ref{klasy} implies
$R_x -$ class is one element provided $\ran (x)\nsubseteq\dom
(a)$. For all such $x$, $\dom (x)$ can be chosen arbitrary. Denote
by $m$ the number of such points $i\in\dom(x)$ that
$x(i)\notin\dom(a)$. Then $1\leq m\leq k$. For a fixed $m$ set
$\ran(x)\backslash\dom(a)$ can be chosen in
$\displaystyle\binom{n-p}{m}$ ways. Similarly, set
\mbox{$\ran(x)\cap\dom(a)$} can be chosen in
$\displaystyle\binom{p}{k-m}$ different ways. Hereby $\ran(x)$ can
be defined in
$\displaystyle{\overset{k}{\underset{m=1}{\sum}}}\binom{n-p}{m}\binom{p}{k-m}$
ways. Then the number of one element $\mathcal{R} -$ classes
equals
\begin{displaymath}
\sum_{k=0}^{n}\sum_{m=1}^{k}\binom{n-p}{m}\binom{p}{k-m}\binom{n}{k}k!\,
\end{displaymath}
where $\displaystyle\binom{n}{k}$ is the number of ways to chose
$\dom (x)$, $k! $ is the number of different $x$ when $\dom (x)$
and $\ran (x)$ are defined.

By Theorem~\ref{klasy} each multi-element $\mathcal{R} -$class is
defined by the domain of its representative $x$, moreover, $|R_x
|>1$ if and only if $\ran(x)\subseteq\dom(a)$. Thus for a fixed
$k$ the number of multi-element $\mathcal{R} -$classes is
$\displaystyle\binom{n}{k}$ and $1\leq k\leq p$. It is clear the
total number of multi-element classes is equal to
$\displaystyle{\overset{p}{\underset{k=1}{\sum}}}\binom{n}{k}$.

Now count the cardinality of class $R_x$. By
Theorem~\ref{klasy}.1) $y\in R_x$ if and only if $\dom (y)=\dom
(x)$ and $\ran (y)\subseteq\dom (a)$. So $\ran(y)$ can be chosen
in $\displaystyle\binom{p}{k}$ ways. Then one defines the map from
$\dom (x)$ to $\ran (y)$ in $k!$ ways. Hence $\displaystyle
|R_x|=\binom{p}{k} k!=[p]_k$.
\end{proof}

\section{Green's relation on $(\mathcal{T}_n ,*_a)$.}
Let $\mathcal{T}_n $ be full symmetric semigroup of all
transformations of the set $\{1,2,\dots,n\}$. For any
$x\in\mathcal{T}_n$ denote by $\ran (x)$ the range of
transformation $x$. The value $|\ran (x)|$ is called the range of
$x$ and is denoted by $\rank (x)$.

Denote by $\rho_x$ the partition of the set $\{1,2,\dots,n\}$
induced by transformation $x$, that is $i$ and $j$ belong to the
same block of the partition $\rho_x$ provided $x(i)=x(j)$. By
$x^{-1} (i) $ denote the full pre-image of the point
$i\in\ran(x)$.

\begin{theorem}
\label{klasyT} Let $n>1$ and $x\in (\mathcal{T}_n ,*_a)$
\begin{enumerate}[1)]
    \item If $\rank (x)\leq \rank (a)$ and for every block $M$
of the partition $\rho_a$, $|\ran (x)\cap M|\leq 1$ then
\setcounter{equation}{8}
        \begin{align*}
    R_x = \{y\:|\:\rho_y=\rho_x \text{ and for every block }M\text{ of the partition } \rho_a,
     \\ |\ran (y)\cap M|\leq
    1\}
    \end{align*}
    and $|R_x|>1$. Otherwise $R_x =\{x\}$.

    \item If $\rank(a)>1,$ $\rank(x)\leq \rank(a)$ and for every block $B_x$ of the partition
    $\rho_x$,
     $B_x\cap\ran (a)\neq \varnothing$ then
     \begin{equation}
    \nonumber
    \begin{split}
    L_x = \{y\:|\ran (y)=\ran (x) \text{ and for every block }B_y\text{ of the partition }
     \rho_y, \\B_y\cap\ran (a)\neq \varnothing \} .
    \end{split}
    \end{equation}

    Otherwise $L_x =\{x\}$.

    If $\rank(a)= 1$ then
    all $\mathcal{L} -$classes are one element.

    \item If $\rank(x)\leq \rank(a)$, for every block $M$ of the partition $\rho_a$,
    \mbox{$|\ran (x)\cap M|\leq 1$}
    and for every block $B_x$ of the partition $\rho_x$, $B_x\cap\ran (a)\neq \varnothing$ then
    \begin{align*}
        H_x = \{y\,|\, \rho_y=\rho_x,\;\ran (y)=\ran (x), \text{ for every block } M \text{ of the partition }
    \rho_a, \\ |\ran (y)\cap M|\leq 1 \text{ and for every block }B_y\text{ of the partition }
    \rho_y,
    \\B_y\cap\ran (a)\neq \varnothing\}.
    \end{align*}
   Otherwise $H_x =\{x\}$.

    \item If $\rank(x)\leq \rank(a)$, for every block $M$ of the partition
    $\rho_a$, $|\ran (x)\cap M|\leq 1$ and either there exists block $B_x$ of the partition
    $\rho_x$, such that $B_x\cap\ran (a)=\varnothing$ or $\rank(a)=1$ then $D_x
    =R_x$;

   if $\rank(x)\leq \rank(a)$, for every block $M$ of the partition $\rho_a$, $|\ran (x)\cap M|>1$
   and for every block $B_x$ of the partition $\rho_x$, $B_x\cap\ran (a)\neq \varnothing$ then $D_x
   =L_x$;

   if $\rank(x)\leq \rank(a)$, for every block $M$ of the partition $\rho_a$, $|\ran (x)\cap M|\leq 1$
   and for every block $B_x$ of the partition $\rho_x$, $B_x\cap\ran (a)\neq \varnothing$ then
    \begin{equation}
    \nonumber
    \begin{split}
    D_x =
    \{y\:|\text{ for every block } B_y \text{ of the partition }
    \rho_y, \:B_y\cap\ran (a)\neq \varnothing,\\ \text{ and for every block } M
\text{ of the partition } \rho_a, \;|\ran (y)\cap M|\leq
    1\};
    \end{split}
    \end{equation}

    in all other cases $D_x =\{x\}$.
\end{enumerate}

\end{theorem}

\begin{proof}
1) Let different elements $x$ and $ y$ from semigroup
$(\mathcal{T}_n,*_a)$
 belong to the same $\mathcal{R} -$ class. By Remark~\ref{umovy} this holds if and only if
 there exist such $u$ and $v$ from $(\mathcal{T}_n ,*_a)$ that:
 \begin{equation}
x=y*_a u, \label{$i'$}
\end{equation}
 \begin{equation}
 y=x*_a v. \label{$i''$}
 \end{equation}

 Since there holds $\rank (xa)\leq\rank
 (a)$ for any transformation $x\in\mathcal{T}_n$ the condition~(\ref{$i'$})
 implies $\rank (x)=\rank(yau)\leq\rank (y)$. Analogously, from~(\ref{$i''$}) one gets
 $\rank (y)=\rank(xav)\leq\rank (x)$. Hence for all $y\in R_x$
 there holds an equality, $\rank (x)=\rank (y)$, moreover
 $\rank(x)\leq\rank(a)$.
This means the points from $\ran (x)$ belong to different blocks
of the partition
 $\rho_a $, that is for every block $M$ of the partition $\rho_a$, $|\ran (x)\cap
M|\leq 1$. So $R_x\subseteq P$, where $P$ is a set from the right
hand side of $(9)$. Now assume that $\rank (x)\leq \rank (a)$ and
for every block $M$ of the partition $\rho_a$ there holds
inequality $|\ran (x)\cap M|\leq 1$. Consider $y\in P$ and $u\in
\mathcal{T}_n $ which is defined on every point $(ya)(i)$ as
$x(i)$ and is defined arbitrary on other points. Analogously,
choose $v\in \mathcal{T}_n $ which is defined on every point
$(xa)(i)$ as $y(i)$. Then the equalities ~(\ref{$i'$}) and
~(\ref{$i''$}) hold. Hereby $y\in R_x$, and the reverse statement
is proved. In this case $|R_x|>1$.

If at least one of the conditions on $x$ fails the above yields
$|R_x| =1$.

2) Let $\rank(a)>1$ and different elements, $x$ and $y$, from
$(\mathcal{T}_n,*_a)$ belong to the same $\mathcal{L} -$ class. By
Remark~\ref{umovy} this holds if and only if there exist such $u$
and $v$ from $(\mathcal{T}_n ,*_a)$ that
\begin{equation}
x=u*_a y \label{$j'$}
\end{equation}
\begin{equation} y=v*_a x
\label{$j''$}
 \end{equation}
Hence $\rank(x)=\rank(y)$. Since $\rank (x)\leq\rank (a)$
equalities~(\ref{$j'$}) and~(\ref{$j''$}) immediately imply $\ran
(x)=\ran (y)$. The last means $\rank (ax)=\rank (x)$ that is for
every block $B_x$ of the partition $\rho_x$, $B_x\cap\ran (a)\neq
\varnothing$. Then for all $y\in L_x$ one gets that for every
block $B_y$ of the partition $\rho_y$, $B_y\cap\ran (a)\neq
\varnothing$. Thus there is an inclusion $L_x\subseteq Q$ where $Q
$ is a set from the right hand side of $(10)$. Conversely, let for
$x$ there hold $\rank(x)\leq \rank(a)$, for every block
 $B_x$ of the partition $\rho_x$, $B_x\cap\ran
(a)\neq \varnothing$, and let $y\in Q$. Then the statement that
for every block $B_y$ of the partition $\rho_y$, $B_y\cap\ran
(a)\neq \varnothing$ implies $\ran (ay)=\ran (y)$. From
$\ran(y)=\ran(x)$ one gets $\ran (ay)=\ran (x)$. So there exist
transformations $u$ and $v$, satisfying the following conditions:
for any $i\in N$ \ $u(i)\in (ay)^{-1} (x(i))$ and respectively for
any $i\in N$ \ $v(i)\in (ax)^{-1} (y(i))$.

The straightforward check shows that $u$ and $v$ satisfy
equalities~(\ref{$j'$}) and ~(\ref{$j''$}). If the conditions of
the statement $2)$ of the theorem fail then $L_x=\{x\}$.

Let $\rank(a)=1$. Assume that $L_x -$ class contains $y\neq x$.
Then ~(\ref{$j'$}) and ~(\ref{$j''$}) imply $\rank (x)=\rank
(y)=1$ and $\ran (x)=\ran (y)$. The last contradicts our
assumption. Hence for all $x$ $|L_x |=1$.

3) Statement about $H_x -$ class follows from 1) and 2) statements
of the theorem and $\mathcal{H} -$ relation definition.

4) As it is known, $x\mathcal{D}y$ if and only if there exist
 $z\in(\mathcal{T}_n,*_a)$ such that $x\mathcal{L}z$ and
 $z\mathcal{R}y$.
Consider possible cases.
\begin{enumerate}[a)]
    \item $|L_x|=1$ and $|R_x|>1$. By statements $1)$ and $2)$ this means that
    $\rank(x)\leq \rank(a)$, for every block $M$ of the partition
    $\rho_a$, $|\ran (x)\cap M|\leq 1$ and either there exists block $B_x$ of the partition
    $\rho_x$,
    such that $B_x\cap\ran (a)= \varnothing$ or $\rank(a)=1$. Then $D_x
    =R_x$.
    \item $|L_x|>1$ and $|R_x|=1$. In this case $\rank(a)>1$, $\rank(x)\leq
    \rank(a)$,
for every block $M$ of the partition $\rho_a$, $|\ran (x)\cap
M|>1$ and for every block $B_x$ of the partition $\rho_x$,
$B_x\cap\ran (a)\neq \varnothing$. Then $D_x =L_x$;
    \item $|L_x|>1$ ³ $|R_x|>1$. Then $\rank(a)>1$, $\rank(x)\leq \rank(a)$, for every block $M$
     of the partition $\rho_a$, $|\ran (x)\cap M|\leq 1$
   and for every block $B_x$ of the partition $\rho_x$, $B_x\cap\ran (a)\neq
   \varnothing$.
   In this case $y\mathcal{D}x$ if and only if for every block $B_y$ of the partition
   $\rho_y$,
$B_y\cap\ran (a)\neq \varnothing$ and for every block $M$ of the
partition $\rho_a$, $|\ran (y)\cap M|\leq
    1$.
    \item If $|L_x|=1$ and $|R_x|=1$ then obviously $D_x =\{x\}$.
\end{enumerate}
\end{proof}

Let $S(n,k)$ be Stirling's number of the second type, that is the
number of (unordered) decompositions of an $n-$element set into
$k$ subsets.
\begin{proposition} Let $\mathcal{T}=(\mathcal{T}_n ,*_a)$,
$n>1$ and $p=\rank (a)$.
\begin{enumerate}[1.]
\item If $p=1$ then all $\mathcal{L} -$classes in $\mathcal{T}$
are one element. The number of one element $\mathcal{L} -$classes
is equal to $n^n$.

Let $p>1.$ Then in $\mathcal{T}$ the number of one element
$\mathcal{L} -$ classes equals
\begin{displaymath}n^n-{\overset{p}{\underset{m=1}{\sum}}}\binom{n}{m}S(p,m)
{\overset{m}{\underset{j=1}{\sum}}}S(n-p,j)\binom{m}{j}j!;
\end{displaymath}
the number of multi-element $\mathcal{L} -$ classes equals
$\displaystyle {\overset{p}{\underset{m=1}{\sum}}}\binom{n}{m}$,
moreover, there are $\displaystyle\binom{n}{m}$ multi-element
$\mathcal{L} -$ classes of the cardinality
\begin{displaymath}S(p,m){\overset{m}{\underset{j=1}{\sum}}}S(n-p,j)\binom{m}{j}j!,\;1\leq
m\leq p.
\end{displaymath}
\item Let $\{a_1,\dots,a_p\} $ be the range of transformation
$a\in \mathcal{T}$.
 Denote $n(a_i)=|a^{-1}(a_i)|$, $1\leq i\leq p$. The number of one element $\mathcal{R}
 -$ classes is equal to
\begin{displaymath}n^n-{\overset{p}{\underset{m=1}{\sum}}}S(n,m){{\underset{1\leq i_1<i_2<\cdots<i_m\leq p}{\sum}}}
 n(a_{i_1})\cdots
n(a_{i_m})m!;
\end{displaymath}
the number of multi-element $\mathcal{R} -$ classes is equal to
$\displaystyle {\overset{p}{\underset{m=1}{\sum}}}S(n,m)$ moreover
there are $S(n,m)$ multi-element $\mathcal{R} -$
 classes of the cardinality
 \begin{displaymath}{{\underset{1\leq i_1 <i_2 <\cdots <i_m\leq
p}{\sum}}}n(a_{i_1})\cdots n(a_{i_m})m!.
\end{displaymath}
\end{enumerate}
\end{proposition}
\begin{proof} 1. The case when $p\leq 1$ is considered in the proof of statement $2)$
of Theorem~\ref{klasyT}.

Let $p>1$. Denote $\rank (x) =m$. By statement $2)$ of the Theorem
~\ref{klasyT} the multi-element $\mathcal{L} -$class in semigroup
$\mathcal{T}$ is uniquely defined by the range of its
representative, $x$, moreover $\rank(x)\leq\rank(a)$, that is
$1\leq m\leq p$. So the number of multi-element  $L -$classes
equals $\displaystyle
{\overset{p}{\underset{m=1}{\sum}}}\binom{n}{m}$. Calculate the
cardinality of this class if $m$ is fixed. By Theorem~\ref{klasyT}
$y\in L_x$ provided $\ran (y)=\ran (x)$ and for every block $B_y$
of the partition $\rho_y$, $B_y\cap\ran a\neq \varnothing$. Hereby
in every block of the partition $\rho_y$ there is at least one
point from $\ran(a)$. The number of distributions of the points
from $\ran (a)$ into the blocks of the partition $\rho_y$, equals
$S(p,m)$. Transformation $y$ maps other $n-p$ points of set
$\{1,2,\dots,n\}$ in $\displaystyle
{\overset{m}{\underset{j=1}{\sum}}}S(n-p,j)\binom{m}{j}j!$ ways,
where $\displaystyle \binom{m}{j}$ is the number of ways to choose
blocks in which $n-p$ points are distributed in $S(n-p,j)$ ways.
Moreover each time in $j!$ ways the chosen blocks can be shuffled.
Hence the total number of transformations in a single
multi-element $\mathcal{L} -$class equals
\begin{displaymath}
S(p,m){\overset{m}{\underset{j=1}{\sum}}}S(n-p,j)\binom{m}{j}j!.
\end{displaymath} Remaining elements from $\mathcal{T}$ form the
set of one element $\mathcal{L} -$classes. Their number is equal
to
\begin{displaymath}
n^n-{\overset{p}{\underset{m=1}{\sum}}}\binom{n}{m}S(p,m)
{\overset{m}{\underset{j=1}{\sum}}}S(n-p,j)\binom{m}{j}j!.
\end{displaymath}

2. By statement $1)$ of Theorem~\ref{klasyT} each multi-element
$\mathcal{R} -$class is defined by the partition $\rho $ of the
set $\{1,2,\dots,n\}$ such that the number of partition blocks is
less or equal to $\rank(a)$. For every $1\leq m\leq p$ there are
$S(n,m)$ decompositions of the set $\{1,2,\dots,n\}$ into $m$
unordered blocks.

Calculate the number of elements in $\mathcal{R} -$class defined
by the partition into blocks $X_1,X_2,\dots,X_m $. By statement
$1)$ of Theorem~\ref{klasyT} this number depends solely on the
number of blocks in the partition, moreover for every $1\leq i\leq
p$ \ $a^{-1}(a_i)$ contains the range of at most one of these
blocks. Let $a_{i_1},a_{i_2},\dots, a_{i_m}$ be such elements from
$\ran(a)$ that $a^{-1}(a_{i_j})$ contains the range of a certain
block. Then there are $m!$ different ways to map the blocks
$X_1,X_2,\ldots,X_m $ to sets
$a^{-1}(a_{i_1}),\ldots,a^{-1}(a_{i_m})$. Then the total number of
the ways to map $X_1,X_2,\dots ,X_m $ to
$a^{-1}(a_{i_1}),\ldots,a^{-1}(a_{i_m})$ is equal to
$n(a_{i_1})\cdots n(a_{i_m})m!$. Hence for a fixed $m$ the
cardinality of a multi-element $\mathcal{R} -$class is equal to
\begin{displaymath}{{\underset{1\leq i_1 <i_2 <\cdots <i_m\leq
p}{\sum}}}n(a_{i_1})\cdots n(a_{i_m})m!.
\end{displaymath}
The total number of elements in all multi-element $\mathcal{R}
-$classes equals
\begin{displaymath}{\overset{p}{\underset{m=1}{\sum}}}S(n,m){{\underset{1\leq
i_1<i_2<\cdots<i_m\leq p}{\sum}}}
 n(a_{i_1})\cdots
n(a_{i_m})m!.
\end{displaymath}
Now it is clear that the number
of one element $\mathcal{R} -$classes is equal to
\begin{displaymath}n^n-{\overset{p}{\underset{m=1}{\sum}}}S(n,m){{\underset{1\leq
i_1<i_2<\cdots<i_m\leq p}{\sum}}}
 n(a_{i_1})\cdots
n(a_{i_m})m!.
\end{displaymath}
\end{proof}

\vspace{0.1cm} \noindent Department of Mechanics and
Mathematics,\\
Kiev Taras Shevchenko University,\\
 64, Volodymyrska st., 01033,
Kiev, Ukraine,\\ e-mail: {\em gtsyaputa\symbol{64}univ.kiev.ua}
\end{document}